\documentclass[11pt]{article}
\usepackage{geometry}                
\usepackage{graphicx}
\usepackage{amssymb}
\usepackage{epstopdf}
\usepackage{amsmath, amsthm}
\usepackage{longtable}
\newtheorem{theorem}{Theorem}
\newtheorem{proposition}[theorem]{Proposition}

\newtheorem{corollary}[theorem]{Corollary}

\newcommand{\Tr}{\mathrm{Tr}}

\begin{document}
\title{The averaged characteristic polynomial for the Gaussian and chiral Gaussian ensembles with a source}
\author{Peter J. Forrester}
\date{}
\maketitle

\noindent
\thanks{\small
Department of Mathematics and Statistics,
The University of Melbourne,
Victoria 3010, Australia email: { P.Forrester@ms.unimelb.edu.au} 
}

\begin{abstract}
\noindent In classical random matrix theory the Gaussian and chiral Gaussian random matrix models with a source are realized as
shifted mean Gaussian, and chiral Gaussian, random matrices with real $(\beta = 1)$, complex ($\beta = 2)$ and real quaternion $(\beta = 4$)
elements. We use the Dyson Brownian motion model to give a meaning  for general $\beta > 0$. In the Gaussian case a
further construction valid for $\beta > 0$ is given, as the eigenvalue PDF of a recursively defined random matrix ensemble. In the case of real or complex
elements, a combinatorial argument is used to compute the averaged characteristic polynomial. The resulting functional forms are shown
to be a special cases of duality formulas due to Desrosiers. New derivations of the general case of Desrosiers' dualities are given.
A soft edge scaling limit of the averaged characteristic polynomial is identified, and an explicit evaluation in terms of so-called
incomplete Airy functions is obtained.
\end{abstract}

\section{Introduction}

It has long been observed that random matrix theory links aspects of probability, combinatorics, integrable systems and asymptotics. As a concrete example, to be further developed in the present paper, let us consider the mean value of the characteristic polynomial of a matrix from the Gaussian orthogonal ensemble (GOE). With $X$ an $N\times N$ matrix of standard Gaussians N$[0,1]$, the latter consists of symmetric matrices $G= \frac{1}{2} (X+X^T)$. It is simple to see \cite{FG06} that the problem of computing $\langle \det (\lambda I_N -G)\rangle_{G\in {\rm GOE}}$ can be reduced to the combinatorial task of counting the number of permutations on $N$ letters consisting of $j$ 2-cycles and $N- 2j$ fixed points, weighted by $(-1)^j \lambda^{N-2j}$ and summed over $j$, and as a consequence
\begin{align}
\label{1} \langle \det (\lambda I_N-G) \rangle_{G\in {\rm GOE}} = 2^{-N} H_N(\lambda ),
\end{align}
where $H_n(x)$ denotes the classical Hermite polynomial. In fact (\ref{1}) holds for any $G=\frac{1}{2}(X+X^T)$ with the elements of $X$ all independent and having mean zero and standard deviation unity.

On the other hand, from the viewpoint of integrable systems, we know that the RHS of (\ref{1}) satisfies a second order constant coefficient differential equation. This differential equation has a turning point at $\lambda= \sqrt{2N}$, and this leads to the asymptotic formula (see e.g.~\cite{Fo93a})
\begin{align}
\label{2} \lim_{N\to\infty} \frac{1}{C_N^{(1)}} e^{-\lambda^2/2} \langle \det(\lambda I_N -G) \rangle_{G\in {\rm GOE}} \Big|_{\lambda= \sqrt{2N} +2^{-1/2}N^{-1/6}y}= 
{\rm Ai}(y),
\end{align}
where Ai$(y)$ denotes the Airy function and
\begin{align}
\nonumber C_N^{(1)} = \pi^{1/4} 2^{-N/2+1/4} (N!)^{1/2} N^{-1/12}.
\end{align}
As $\sqrt{2N}$ is the leading boundary of support for the eigenvalues, this scaling of the characteristic polynomial is thus seen to provide a statistical probe of the largest eigenvalues. It is furthermore true that the RHS of (\ref{1}) permits an integral form, which implies the so-called duality formula \cite{BF97a}, \cite{BH01}
\begin{align}\label{3}
\langle \det (\lambda I_N- G)\rangle_{G\in {\rm GOE}} = \langle (\lambda+ ix)^N \rangle_{x\in {\rm N}[0,1/\sqrt{2}]},
\end{align}
and this allows for an alternative derivation of (\ref{2}).

Our aim in this paper is to seek analogous inter-relations for so called matrix models with a source. Both Gaussian and chiral Gaussian classes of such matrices will be considered. In the Gaussian case, a simple example of such a matrix is $G+ G_0$, where $G\in$ GOE and $G_0$ is the diagonal matrix $\mathrm{diag}\: [(\mu)^r, (0)^{N-r}]$, $(a)^p$ denoting $a$ repeated $p$ times. The matrix $G_0$ is referred to as a rank $r$ spiked matrix. Similarly, a simple example of a
matrix model with a source in the 
chiral Gaussian class is
\begin{align}\label{K}
 \left[ \begin{array}{cc}
\mathbf{0}_{n \times n} & X_{n \times p}\\
(X^T)_{p \times n} & \mathbf{0}_{p \times p}
\end{array}\right] + \left[ \begin{array}{cc}
\mathbf{0}_{n \times n} & (X_0)_{n \times p}\\
(X_0^T)_{p \times n} & \mathbf{0}_{p \times p}
\end{array}\right] =: Z + Z_0,
\end{align}
where the rectangular matrix $X_{n \times p}$ consists of standard Gaussians and the entries of $(X_0)_{n \times p}$ are equal to $\mu$ on
the diagonal for the first $r$ rows, and zero elsewhere.

We begin in Section 2 by using combinatorial arguments to derive the generalization of (\ref{1}) for both the Gaussian and chiral Gaussian ensembles with
a source, in the case the elements are real (labelled $\beta = 1$), and the case the elements are complex (labelled $\beta = 2$).
In Section 3 we use the Dyson Brownian motion model \cite{Dy62b} to give meaning to the Gaussian ensemble with a source for
general parameter $\beta > 0$, and in Section 4 we do similar for the chiral Gaussian ensemble with a source. In Section 5 duality formulas from
\cite{De08} are introduced, which generalize (\ref{3}) to the case of the Gaussian and chiral Gaussian ensembles with a source. The main purpose of this
section is to give new derivations of the duality formulas. In Section 6 we link the evaluations of the averaged characteristic polynomials from
Section 2 with (a special case of) the duality formulas from Section 5, and we proceed to compute soft edge scaled limits, giving rise to a generalization of (\ref{2}).

\section{The averaged characteristic polynomial for Gaussian and chiral Gaussian ensembles with
a source: real or complex entries}
\setcounter{equation}{0}
\subsection{The shifted mean GOE and GUE}\label{s2.1}
Let $H$ be a symmetric matrix, and require that all entries on and above the diagonal be independent with mean zero, and require that the
variance of the elements above the diagonal equal $1/2$. Let $H_0$ be a fixed diagonal matrix, $H_0 = {\rm diag} \,(s_1,s_2,\dots,s_n) $.
We seek the mean value of the characteristic polynomial of the matrix $H+H_0$.

\begin{proposition}
Let $H$ and $H_0$ be as specified above. We have
\begin{equation}\label{pp1}
\langle \det (\lambda I_N - (H + H_0)) \rangle =
\Big \langle \prod_{j=1}^N ( \lambda -s_j + i x ) \Big \rangle_{x \in {\rm N}[0,1/\sqrt{2}]}.
\end{equation}
\end{proposition}

\noindent
{\it Proof. } 
Let $H = [h_{ij}]_{i,j=1,\dots,N}$. We have
\begin{equation}\label{6.4}
\det (\lambda I_N - (H + H_0)) =
\sum_{P \in S_N} \varepsilon(P) \prod_{l=1}^N ( \lambda_{l,P(l)} - s_{l,P(l)}  - h_{l,P(l)}),
\end{equation}
where $\varepsilon(P)$ denotes the parity of $P$ and
$$
\lambda_{i,j} - s_{i,j} = \left \{
\begin{array}{ll} \lambda - s_i, & i = j \\
0, & i \ne j. \end{array}
\right.
$$
The fact the entries of $H$ each have mean zero tells us that after averaging the only terms in (\ref{6.4}) to give a non-zero contribution consist
entirely of fixed points and 2-cycles. Each fixed point $P(l) = l$ contibutes $\lambda - s_l$, while each 2-cycle contributes $1/2$. Let there be
$N - 2j$ fixed points, and $j$ 2-cycles. The (weighted) number of ways to form $j$ 2-cycles from $2j$ numbers is
$$
2^{-j} {(2j)! \over 2^j j!},
$$
while the (weighted) number of ways to choose the fixed points is
$$
\sum_{1 \le l_1 < l_2 < \cdots < l_{N-2j} \le N}
\prod_{k=1}^{N-2j} ( \lambda - s_k).
$$
Further, such permutations have parity $(-1)^j$ and so
\begin{equation}\label{pp}
\langle \det ((\lambda I_N - (H + H_0)) \rangle =
\sum_{j=0}^{[N/2]} (-1)^j  2^{-j} {(2j)! \over 2^j j!}
\sum_{1 \le l_1 < l_2 < \cdots < l_{N-2j} \le N}
\prod_{k=1}^{N-2j} (\lambda - s_{l_k}).
\end{equation}
Using the facts that
$$
{1 \over \sqrt{\pi} } \int_{-\infty}^\infty e^{-x^2} x^{2j} \, dx = 2^{-j} {(2j)!  \over 2^j j!}, \qquad \int_{-\infty}^\infty e^{-x^2} x^{2j+1} \, dx =0,
$$
we see that (\ref{pp}) can equivalently be written in the form (\ref{pp1}). \hfill $\square$

\medskip

In the special case that $s_j = 0$ for each $j$, we see that (\ref{pp1}) reduces to (\ref{3}).

Similar reasoning applies to the calculation of the average of the characteristic polynomial for the shifted
mean Gaussian unitary ensemble (GUE). In this setting $H$ is required to be an Hermitian matrix with all entries on and above the
diagonal independent with mean zero, and with the variance of the sum of the real and imaginary parts of the off diagonal
equalling $1/2$.  The averaged characteristic polynomial is again given by (\ref{pp1}).

\subsection{The shifted mean chGOE and chGUE}
Consider the matrix structure (\ref{K}) and suppose the entries of $X$ are all independent with mean zero and variance unity. 
Let $X_0$ have entries on the diagonal equal to $s_j$ for rows $j=1,\dots,p$, and all other entries equal to zero.    
In the case that $s_j = 0$ for each $j$ and thus $Z_0 = {\bf 0}_{(n+p) \times (n + p)}$ in (\ref{K}) we know that \cite{FG06}
\begin{equation}\label{fg}
\langle \det (\lambda I_{n+p} - Z ) \rangle =    (-1)^p p! \lambda^{n-p} L_p^{n-p}(\lambda^2),
\end{equation}                                                        
where $L_m^a(x)$ denotes the classical Laguerre polynomial. We seek to generalize this formula to the case of nonzero $s_j$.

\begin{proposition}
Let $X$ and $X_0$ be as specified above, and let $Z$ and $Z_0$ then be specified by (\ref{K}). We have
\begin{equation}\label{pp1A}
\langle \det (\lambda I_{n+p} - (Z + Z_0)) \rangle =
{(-1)^p \lambda^{n-p} e^{\lambda^2} \over \Gamma(n-p+1)}
\int_0^\infty t^{n-p} e^{-t} \, {}_0 F_1(n-p+1;-\lambda^2 t)
\prod_{l=1}^p (t +  s_l^2) \, dt.
\end{equation}
\end{proposition}

\noindent
{\it Proof. } 
In the analogue of (\ref{6.4}) for the matrix structure (\ref{K}) it is again true that after averaging the only permutations giving a non-zero
contribution are those consisting entirely of fixed points and 2-cycles. Further, there must be at least $n-p$ fixed points, due to (\ref{K}) having
$n-p$ zero eigenvalues.

The 2-cycles $(j_1 j_2)$ giving a non-zero contribution have $j_1$ ($j_2$) corresponding to a row (column) of $X$. Each time
$(j_1 j_2)$ is a diagonal entry $x_{j_1 \,p+j_1} - s_{j_1}$ there is a choice in the analogue of (\ref{6.4}) of averaging $s_{j_1}^2$
(contribution $s_{j_1}^2$) or of averaging $x_{j_1 \, j_2}^2$ (contribution unity). In particular, if there are $j$ 2-cycles (and thus
$n+p - 2j$ fixed points) we must catalogue the 2-cycles according to the occurrences of there being zero, one, two etc.~diagonal
entries each with the choice of averaging $s_j^2$ for the corresponding $j$. Doing this we can see that the weighted contribution to the
averaged characteristic polynomial is
\begin{eqnarray}\label{av}
&& (-1)^j \lambda^{n+p - 2j} \Big \{ 
\binom{n}{j}  \binom{p}{j}  j! +
\Big ( \sum_{l=1}^p s_l^2 \Big ) \binom{n-1}{j-1}  \binom{p-1}{j-1}  (j-1)! \nonumber \\
&& \quad + \Big ( \sum_{1 \le l_1 < l_2 \le p} s_{l_1}^2 s_{l_2}^2 \Big )  \binom{n-2}{j-2}  \binom{p-2}{j-2}  (j-2)! + \cdots +
 \sum_{1 \le l_1 < \cdots  < l_p \le p} s_{l_1}^2 \cdots s_{l_p}^2 \Big \}
\end{eqnarray}

The above consists of $p+1$ terms. Labelling them $r=0,1,\dots,p$, the corresponding combinatorial factor
$$
  \binom{n-r}{j-r}  \binom{p-r}{j-r}  (j-r)! 
$$
results from counting the number of ways of choosing $j-r$ rows from $n-r$, independently choosing $j-r$ columns from $p-r$, and
further independently ordering the $j-r$ rows.

The averaged characteristic polynomial comes from summing (\ref{av}) from $j=0$ to $p$. For given $j$, we replace $j \mapsto p - j$.
Making use of the explicit power series form of the
Laguerre polynomial,
$$
L_p^a(x) = \sum_{j=0}^p  \binom{p+a}{p-j}  {(-1)^j \over j!} x^j,
$$
we read off that
\begin{eqnarray}\label{17.1x}
&& \langle \det (\lambda I_{n+p} - (Z + Z_0)) \rangle =
(-1)^p \lambda^{n-p} \Big \{
p! L_p^{n-p}(\lambda^2)  + \Big ( \sum_{l=1}^p s_l^2 \Big ) (p-1)!   L_{p-1}^{n-p}(\lambda^2)   \nonumber \\
&& \qquad + 
 \Big ( \sum_{1 \le l_1 < l_2 \le p} s_{l_1}^2 s_{l_2}^2 \Big )  (p-2)!    L_{p-2}^{n-p}(\lambda^2) +  \cdots
 +  \sum_{1 \le l_1 < \cdots  < l_p \le p} s_{l_1}^2 \cdots s_{l_p}^2 \Big \}.
 \end{eqnarray}
 
 To proceed further, we substitute for the Laguerre polynomials their integral representation
 \begin{equation}\label{LP}
 m! L_m^a(x) = {e^{x} \over \Gamma(a+1)}
 \int_0^\infty e^{-t} t^{m+a}  {}_0 F_1(a+1;- x t) \, dt
 \end{equation}
 (a special case of (\ref{23.1}) below). This allows the RHS of (\ref{17.1x}) to be rewritten
 \begin{eqnarray*}
&& 
(-1)^p {\lambda^{n-p} e^{\lambda^2} \over \Gamma(n-p+1) }
\int_0^\infty t^{n-p} e^{-t}  {}_0 F_1(n-p+1;- \lambda^2 t)  \Big \{ t^p +  t^{p-1}  \Big ( \sum_{l=1}^p s_l^2 \Big )\nonumber \\
&& \qquad +
 t^{p-2}  \Big ( \sum_{1 \le l_1 < l_2 \le p} s_{l_1}^2 s_{l_2}^2 \Big ) + \cdots + 
 \sum_{1 \le l_1 < \cdots  < l_p \le p} s_{l_1}^2 \cdots s_{l_p}^2 \Big \} \, dt.
 \end{eqnarray*}
 Performing the sum gives (\ref{pp1A}). \hfill $\square$
 
 \medskip
 Note that setting $s_l = 0$, $l=1,\dots,p$ in (\ref{17.1x}) reclaims (\ref{fg}).
 
 If we replace $X^T$ in (\ref{K}) by $X^\dagger$ , then we could consider the case of $X$ having complex entries, with real and
 imaginary parts having mean zero and the sum of the variance of the real and imaginary parts equalling unity. The above working
 remains valid and thus the average characteristic polynomial is again given by (\ref{pp1A}).
 
 Generally, for the matrix structure (\ref{K}) with $n \ge p$ there are $n-p$ zero eigenvalues, and the eigenvalues come in $\pm$ pairs.
 Squaring the nonzero eigenvalues gives the eigenvalues of $(X + X_0)^\dagger(X + X_0)$. We can therefore read off from
(\ref{pp1A}) the corresponding averaged characteristic polynomial.

\begin{corollary}\label{C3}
Let the entries of $X$ have mean zero and expected value for their absolute value squared equal to unity. Let $X_0$ be nonzero
down it diagonal only, with $X_0^T X_0 = {\rm diag} \,(\mu_1,\dots,\mu_p)$.
We have
$$
\langle \det (\lambda I_{p} - (X + X_0)^\dagger (X + X_0)) \rangle =
{(-1)^p e^{\lambda} \over \Gamma(n-p+1)}
\int_0^\infty t^{n-p} e^{-t} \, {}_0 F_1(n-p+1;-\lambda t)
\prod_{l=1}^p (t  +  \mu_l) \, dt.
$$
\end{corollary}

\section{Gaussian matrix models with a source}
\setcounter{equation}{0}
\subsection{The three classical cases}

Let $G$ be a random Hermitian matrix with real ($\beta=1$), complex ($\beta=2$) or real quaternion ($\beta=4$) entries. With the joint distribution of the elements of $G$ proportional to $e^{-(\beta/2) \Tr \: G^2}$, such random matrices are said to form the Gaussian orthogonal ($\beta=1$), unitary ($\beta=2$) and symplectic ($\beta=4$) ensembles respectively. Note that in the cases $\beta=1$ and $2$ this definition coincides with the definition $G=\frac{1}{2} (X+X^T)$ for $X$ consisting of standard real ($\beta=1$) or complex ($\beta=2$) entries. By definition the corresponding Gaussian matrix models with a source \cite{BH01} are then specified as random matrices $H= G+M$, where $M$ is Hermitian and has entries from the same number field as $G$. Note that the set of matrices $H$ may equally as well be referred to as the shifted mean Gaussian ensembles.

For future purposes it is convenient to replace $M$ by $e^{-\tau} M$ ($\tau>0$), and to modify the variance in the joint distribution of the elements of $G$ so that the latter is proportional to $e^{-(\beta/2) \Tr\: G^2/(1-e^{-2\tau})}$. We then have that the joint distribution of the elements of $H$ is proportional to
\begin{align}
\label{8.0} e^{-(\beta/2) \Tr\: (H-e^{-\tau}M)^2/(1-e^{-2\tau})}.
\end{align}
We now take up the problem of determining the eigenvalue PDF of $H$ \cite[\S11.1 \& \S 13.4.3]{Fo10}.

A fundamental result in relation to Hermitian matrices with real, complex or quaternion real elements is that the measure for the change of variables corresponding to the diagonalization $H=UDU^\dagger$, where $D$ denotes the diagonal matrix formed by the eigenvalues $\vec{\lambda} := \{ \lambda_j \}_{j=1,...,N}$ and $U$ the matrix of eigenvectors, is proportional to (see e.g. \cite[Prop 1.3.4]{Fo10})
\begin{align}
\label{8.1} \prod_{1\leq j <k\leq N} |\lambda_k - \lambda_j|^{\beta} \prod_{j=1}^{N} d\lambda_j \; (U^{\dag} dU).
\end{align}
Here $(U^{\dag} dU)$ denotes the normalized Haar measure on real orthogonal ($\beta=1$), unitary ($\beta=2$) and unitary symplectic matrices ($\beta=4$). Using (\ref{8.1}) we can then deduce the well known result that the eigenvalue PDF, $G_{\tau} (\vec{\lambda}; \vec{\mu})$ say, with $\vec{\mu} = \{ \mu_j\}_{j=1,...,N}$ denoting the eigenvalues of $M$, for the shifted mean Gaussian ensembles is proportional to
\begin{align}
\label{8.2} \prod_{1\leq j < k\leq N} |\lambda_k-\lambda_j|^{\beta} e^{-\tilde{\beta}\sum_{j=1}^N \lambda_j^2- \tilde{\beta} t^2 \sum_{j=1}^N \mu_j^2} \int e^{2\tilde{\beta} t \Tr\: (U\Lambda U^{\dag}L)} \: (U^{\dag} dU).
\end{align}
Here
\begin{align}
\label{ll} \Lambda = \mathrm{diag} (\lambda_1,...,\lambda_N), \qquad L=\mathrm{diag} (\mu_1,...,\mu_N)
\end{align}
and
\begin{align}
\nonumber \tilde{\beta} = \beta/(2(1-e^{-2\tau})), \qquad  t=e^{-\tau}.
\end{align}

Denoting
\begin{align}
\label{AF} \int e^{\Tr \: (U\Lambda U^{\dag} L)} (U^{\dag} dU)=: {}^{\phantom{/}}_0 {\mathcal F}_0^{2/\beta} (\vec{\lambda}; \vec{\mu})
\end{align}
the (un-normalized) PDF (\ref{8.2}) reads
\begin{equation}
\label{AF1} 
\prod_{1 \le j < k \le N} | \lambda_k -  \lambda_j |^\beta
e^{- \tilde{\beta} \sum_{j=1}^N \lambda_j^2 - \tilde{\beta} t^2  \sum_{j=1}^N \mu_j^2 }
 {}^{\phantom{/}}_0 {\mathcal F}_0^{(2/\beta)} (\vec{\lambda};2  \tilde{\beta}  t  \vec{\mu}).
 \end{equation}
 As it stands (\ref{AF1}) only has meaning for $\beta = 1,2$ and 4. One concern of contemporary random matrix theory is to give
 constructions/ interpretations of such PDFs for general $\beta > 0$ (see e.g.~\cite{Fo11y}). In fact two approaches are possible ---
 one as the solution of a Fokker-Planck equation, which we take up in the next section, and the other as the eigenvalue PDF of an
 ensemble of matrices defined recursively, which is to be considered in the Appendix. We remark that in the case of the general
 variance Wishart ensemble, the corresponding interpretations of the analogue of (\ref{AF1}), proposed as a candidate for 
 the corresponding $\beta$-generalized distribution in \cite{Wa11}, have
 recently been given in \cite{Fo11y}.
 
 \subsection{Brownian evolution interpretation}\label{s2.2}
 
 We see from (\ref{8.0}) that as $\tau \to \infty$ the element PDF of $H$ reduces to that for matrices $G$. The
 corresponding eigenvalue PDF is then well known (see e.g.~\cite[Prop.~1.10.4]{Fo10}) to be proportional to $e^{-\beta U}$ with
 \begin{equation}
\label{U}
U = {1 \over 2} \sum_{j=1}^N \lambda_j^2 - \sum_{1 \le j < k \le N} \log |\lambda_k - \lambda_j |.
\end{equation}
A key observation, due to Dyson \cite{Dy62}, is that $e^{-\beta U}$ can be interpretated as the Boltzmann factor for a classical gas of
$N$ particles on the real line with potential (\ref{U}), interacting at inverse temperature $\beta$. Generally the Brownian evolution of
such a classical gas is specified by the Fokker-Planck equation
\begin{equation}\label{FP}
{\partial p_\tau \over \partial \tau} = {\mathcal L} p_\tau, \qquad
\mathcal L = \sum_{j=1}^N {\partial \over \partial \lambda_j} \Big (
{\partial U \over \partial \lambda_j} + {1 \over \beta} {\partial \over \partial  \lambda_j} \Big ).
\end{equation}

For $\beta = 1,2$ and 4 this equation was shown by Dyson \cite{Dy62} to specify the eigenvalue PDF for Hermitian random
matrices with real ($\beta = 1$), complex ($\beta = 2$) and real quaternion  ($\beta = 4$)  elements, and independent real
and imaginary parts themselves Brownian paths confined by a harmonic potential.
Later (see e.g.~\cite{Ha90}) this specification was shown to be equivalent to choosing the elements of $H$ according to the
distribution (\ref{8.0}). Note that the solution of (\ref{FP}) corresponding to (\ref{8.0}) has the delta function initial condition
$\prod_{l=1}^N \delta(\lambda_l - \mu_l)$. Consequently, for general $\beta > 0$ (\ref{AF1}) is then, up to a
proportionality depending on $\tau$, the Green function solution of the Fokker-Planck equation (\ref{FP}).

Development of the theory of the eigenfunctions and eigenvalues of (\ref{FP}) allows ${}^{\phantom{/}}_0 {\mathcal F}_0^{(2/\beta)} $
in (\ref{AF1}) to be made more explicit. This requires introducing the symmetric Jack polynomials $P_\kappa(z;\alpha)$ where
$z = (z_1,\dots,z_N)$, $\kappa$ denotes a partition of length less than or equal to $N$ (we write $\ell(\kappa) \le N$), and $\alpha$
is a parameter. These polynomials are homogeneous, $P_\kappa(cz;\alpha) = c^{|\kappa|} P_\kappa(z;\alpha)$ for $c$ a scalar,
and can be uniquely specified in terms of the polynomial eigenfunctions of a certain differential operator (see e.g.~\cite[\S 12.6]{Fo10}
for more details). For $\alpha = 1$ the Jack polynomials coincide with the well known Schur polynomials, and for $\alpha = 2$ ($\alpha = 1/2$)
they are the zonal polynomials for the real orthogonal (unitary symplectic) group. Introducing too the quantity $d_\kappa'$ as in
\cite[eq.~(12.60)]{Fo10}, it is a known result that ${}^{\phantom{/}}_0 {\mathcal F}_0^{(2/\beta)} $ in (\ref{AF1}) can be expanded
\begin{equation}
\label{0F0}
{}^{\phantom{/}}_0 {\mathcal F}_0^{(\alpha)}(\vec{x};\vec{y}) = \sum_\kappa \alpha^{|\kappa|}
{P_\kappa(x;\alpha) P_\kappa(y;\alpha) \over d_\kappa' P_\kappa((1)^N;\alpha) },
\end{equation}
making (\ref{AF1}) explicit for all $\beta > 0$.

\section{Chiral Gaussian random matrix models with a source}
\setcounter{equation}{0}
The chiral Gaussian ensembles consist of the first of the matrices in (\ref{K}), with $X^T$ replaced by $X^\dagger$, and the
entries of the $n \times p$ ($n \ge p$) matrix $X$ independent standard Gaussian real $(\beta = 1)$, complex $(\beta = 2)$ and
real quaternion $(\beta = 4)$ random variables. The shifted mean version of this ensemble is then given by the sum of the
matrices (\ref{K}).

With $M := X_0$, as in (\ref{8.0}) it is convenient to introduce a parameter $\tau$, $\tau > 0$, such that $Y = X + M$ is modified to have
joint distribution of its elements proportional to
\begin{equation}
\label{yb}
e^{-(\beta/2) {\rm Tr} \, (Y - e^{-\tau} M)^\dagger (Y - e^{-\tau} M)/(1 - e^{-2 \tau}) }.
\end{equation}
We seek from this the eigenvalue PDF of $(Y^\dagger Y)^{1/2}$. Let the eigenvalues of the latter be denoted
$\{y_j\}_{j=1,\dots,p}$. The strategy (see \cite[\S 11.2.2]{Fo10}) is to introduce the singular value decomposition $Y = U \Lambda V$
where $U$ is an $n \times p$ unitary matrix, $V$ is a $p \times p$ unitary matrix (both $U$ and $V$ have elements of the same
type as $Y$) and $\Lambda$ is an $n \times p$ matrix with diagonal entries $y_1,\dots,y_p$ and all other entries zero. The
measure for this change of variables in proportional to 
$$
\prod_{j=1}^p y_j^{\beta a + 1} \prod_{1 \le j < k \le p}
|y_k^2 - y_j^2|^\beta \prod_{j=1}^p dy_j \, (U^\dagger dU) (V^\dagger dV) ,
$$
where $a=n-p+1-2/\beta$.

Multiplying this and (\ref{yb}), then integrating out over $U$ and $V$ shows that the PDF for $\{y_j\}_{j=1,\dots,p}$,
$G_\tau^{\rm ch}(\vec{y};\vec{\mu})$ say, is proportional to \cite[eq.~(11.105)]{Fo10}
\begin{align}\label{17.1}
\prod_{j=1}^p y_j^{\beta a + 1}& \prod_{1 \le j < k \le p}
|y_k^2 - y_j^2|^\beta e^{-\tilde{\beta} \sum_{j=1}^p y_j^2 - \tilde{\beta} t^2 \sum_{j=1}^p \mu_j^2} \nonumber \\
& \times {}^{\phantom{/}}_0 {\mathcal F}_1^{(2/\beta)}(\beta n / 2; \{ y_j^2 \}_{j=1,\dots,p}; \tilde{\beta} t \{ \mu_j^2 \}_{j=1,\dots,p} )
\end{align}
where $\tilde{\beta}$, $t$ are as in (\ref{8.2}), $\{ \mu_j^2\}$ are the eigenvalues of $(M^\dagger M)^{1/2}$ and
\begin{align}\label{17.1a}
 {}^{\phantom{/}}_0 {\mathcal F}_1^{(2/\beta)} & (\beta n / 2; \{ y_j^2 \}_{j=1,\dots,p}; \tilde{\beta} t \{ \mu_j^2 \}_{j=1,\dots,p} ) \nonumber \\
 & = \int (U^\dagger d U) \int  (V^\dagger d V) \,
 e^{\tilde{\beta} t {\rm Tr} (V \Lambda^\dagger U^\dagger L + L^\dagger U \Lambda V^\dagger)}.
 \end{align}
 If we introduce $x_j = y_j^2$, $m_j = \mu_j^2$ $(j=1,\dots,p)$ as the eigenvalues of $X^\dagger X$ and $M^\dagger M$ respectively, then
 as a PDF in $\vec{x} := \{x_j \}_{j=1,\dots,p}$ (\ref{17.1}) reads
 \begin{align}\label{17.1b}
\prod_{j=1}^p x_j^{\beta a/2} & \prod_{1 \le j < k \le p}
|x_k - x_j|^\beta e^{-\tilde{\beta} \sum_{j=1}^p x_j - \tilde{\beta} t^2 \sum_{j=1}^p m_j} 
{}^{\phantom{/}}_0 {\mathcal F}_1^{(2/\beta)}(\beta a / 2 +\beta(p-1)/2 + 1; \vec{x}_j ; \tilde{\beta} t \vec{m} ).
\end{align}

Analogous to the discussion in \S \ref{s2.2}, the fact that (\ref{yb}) satisfies a multi-dimensional heat equation can be used 
\cite{Fo98e} to characterize (\ref{17.1}) as a solution of the Fokker-Planck equation (\ref{FP}), with $N \mapsto p$,
$\{ \lambda_j \} \mapsto \{ y_j \}$ and
$$
U = {1 \over 2} \sum_{j=1}^p y_j^2 - {a + 1/\beta \over 2} \sum_{j=1}^p \log y_j^2 -
\sum_{1 \le j < k \le p} \log | y_k^2 - y_j^2 |,
$$
defined on the half line $y > 0$. Moreover, by developing the theory of the corresponding eigenfunctions and eigenvalues
\cite{BF97a}, we can show that for general $\alpha > 0$
$$
 {}^{\phantom{/}}_0 {\mathcal F}_1^{(\alpha)}(c;\vec{x};\vec{y}) =
 \sum_\kappa {\alpha^{|\kappa|} P_\kappa(x;\alpha)  P_\kappa(y;\alpha)  \over
 [c]_\kappa^{(\alpha)} d_\kappa' P_\kappa((1)^N;\alpha)}
 $$
 (cf.~(\ref{0F0})) where
 $$
 [c]_\kappa^{(\alpha)} := \prod_{l=1}^N {\Gamma(c + (j-1)/\alpha + \kappa_j) \over
 \Gamma(c + (j-1)/\alpha ) }.
 $$
 Thus (\ref{17.1b}) is explicit for all $\beta > 0$.
 
 \section{Duality formulas}
 \setcounter{equation}{0}
 The duality formula (\ref{3}) for the average of the characteristic polynomial in the GOE can be
 generalized to a duality formula for the $n$-th moment of the characteristic polynomial in the
 GOE \cite{BF97a}, \cite{BH01}. In fact such formulas hold not only for the GOE, but their generalization known as
 the Gaussian $\beta$-ensemble. To define the latter, let ME${}_{\beta,N}(w(x))$ refer to the matrix
 ensemble with eigenvalue PDF proportional to 
 \begin{equation}\label{w1}
 \prod_{l=1}^N w(x_l) \prod_{1 \le j < k \le N} | x_k - x_j |^\beta.
 \end{equation}
 The Gaussian $\beta$-ensemble is then ME${}_{\beta,N}(e^{-x^2})$.
 
 In terms of this notation, we have as a generalization of (\ref{3}) the duality formula
 \cite[eq.~(13.162)]{Fo10}
 \begin{equation}\label{w2}
 \Big \langle \prod_{j=1}^N (x - \sqrt{2/\beta} \lambda_j )^n
 \Big \rangle_{{\rm ME}_{\beta,N}(e^{-\lambda^2})} =
 \Big \langle \prod_{j=1}^n (x - i  \lambda_j )^N
 \Big \rangle_{{\rm ME}_{4/\beta,N}(e^{-\lambda^2})}. 
 \end{equation}
 Note how the role of $N$ and $n$ is interchanged on the two sides, and that on the RHS $\beta \mapsto 4/\beta$.
 Desrosiers \cite{De08} (see also \cite{BH01} in the case $\beta = 2$) has generalized (\ref{w2}) to the case of the
 Gaussian matrix model with a source. These formulas, together with derivations distinct from those given in
 \cite{De08}, will be given in the next sections.
 
 \subsection{Gaussian matrix model with a source}
 Let us denote the eigenvalue PDF (\ref{AF1}) for the Gaussian matrix model with a source in the case
 $\tilde{\beta} = t = 1$, by ME${}_{\beta,N}(e^{-\lambda^2};\vec{\mu})$. In terms of this notation, the duality formula
 of \cite{De08} reads
 \begin{align}\label{Dr1}
  \Big \langle \prod_{j=1}^n \prod_{k=1}^N (s_j - \sqrt{\alpha} y_k)  \Big \rangle_{{\rm ME}_{2/\alpha,N}(e^{-y^2};-i \vec{x})}
   = i^{nN} 
  \Big \langle \prod_{j=1}^N \prod_{k=1}^n (y_k + \sqrt{\alpha} x_j) \Big  \rangle_{{\rm ME}_{2 \alpha,n}(e^{-y^2};-i \vec{s})}.
  \end{align}
   In writing this from the form given in \cite{De08} we have made use of the fact that the unspecified normalization in
  (\ref{AF1}) is independent of $\vec{\mu}$.
  
  The derivation of (\ref{Dr1}) given in \cite{De08} relies on special properties of the so called Dunkl transform.
 Here we will verify (\ref{Dr1}) by using the generalized Hermite polynomials $P_\kappa^{(H)}(\vec{x};\alpha)$ \cite{BF97a}.
  The latter can be specified by the requirements that their expansion in terms of Jack polynomials has the structure
  $$
  P_\kappa^{(H)}(\vec{x};\alpha) = P_\kappa(\vec{x};\alpha) + \sum_{\mu: |\mu| < |\kappa|}
  a_{\kappa \mu} P_\mu(\vec{x};\alpha)
  $$
  for some coefficients $a_{\kappa \mu}$, and that they exhibit the orthogonality
  $$
    \Big \langle P_\kappa^{(H)}(\vec{x};\alpha) P_\rho^{(H)}(\vec{x};\alpha)  \Big \rangle_{{\rm ME}_{2/ \alpha,N}(e^{-x^2})}
    \propto \delta_{\kappa,\rho}.
 $$
 
 We require a number of special properties of the generalized Hermite polynomials, in particular integration and expansion
 formulas. The integration formulas required are \cite{BF97a}
 \begin{align}
 \label{20.1}
 \langle P_\kappa(i \vec{y};\alpha) \rangle_{{\rm ME}_{2/\alpha,N}(e^{-y^2};-i \vec{x})} & =  P_\kappa^{(H)}(\vec{x};\alpha) \\
 \label{20.2}
  \langle P_\kappa^{(H)}( \vec{y};\alpha) \rangle_{{\rm ME}_{2/\alpha,N}(e^{-y^2}; \vec{x})} & =  P_\kappa(\vec{x};\alpha).
  \end{align}
  Similarly two expansion formulas are required.
  The first is  \cite[eq.~(13.130) upon the limiting procedure of Exercises 13.3 Q.5]{Fo10} 
  \begin{equation}\label{20}
  \prod_{k=1}^N \prod_{l=1}^n (x_k - y_l) =
  \sum_{\kappa} (-1)^{|\mu|} P_\kappa^{(H)}(\vec{x};\alpha)
  \alpha^{-|\mu|/2} P_\mu^{(H)}(\sqrt{\alpha} \vec{y};1/\alpha)
  \end{equation}
  where $\mu := (N)^n - \kappa'$, with $\kappa'$ the conjugate of $\kappa$, while the second is its Jack polynomial analogue
  (see e.g.~\cite[eq.~(12.212)]{Fo10})
  \begin{equation}\label{20a}
  \prod_{k=1}^N \prod_{l=1}^n (x_k + y_l) = \sum_\kappa P_\kappa(x;\alpha) P_\mu(y;1/\alpha).
  \end{equation}
  
  Consider (\ref{20}) rewritten to read
  $$
  \prod_{j=1}^n \prod_{k=1}^N (s_j - \sqrt{\alpha} y_k) =
  \sum_\kappa (-1)^{|\kappa|}
  (\sqrt{\alpha})^{|\kappa|}   P_\kappa^{(H)}(\vec{y};\alpha)  P_\mu^{(H)}(\vec{s};\alpha).
  $$
  Substituting in the LHS of (\ref{Dr1}) and making use too of (\ref{20.2}) shows that the LHS of (\ref{Dr1}) can be written
  \begin{equation}\label{c1}
  \sum_\kappa i^{|\kappa|} (\sqrt{\alpha})^{|\kappa|}  P_\kappa(\vec{x};\alpha)  P_\mu^{(H)}(\vec{s};\alpha).
    \end{equation}
 Next we consider (\ref{20a}) rewritten to read
 $$
 \prod_{j=1}^N \prod_{k=1}^n ( y_k + \sqrt{\alpha} x_j) =
 \sum_\kappa P_\kappa(\sqrt{\alpha} \vec{x};\alpha)
 P_\mu(i \vec{y};1/\alpha) (-i)^{|\mu|}.
 $$
 Substituting in the RHS of (\ref{Dr1}) and making use of (\ref{20.1}) shows that the former
 is also given by (\ref{c1}), thus establishing the validity of (\ref{Dr1}).
 
 \subsection{Chiral Gaussian matrix model with a source}
 With $\tilde{\beta} = t = 1$, $\beta a/2  \mapsto a $, we will denote the transformed PDF for
 the chiral Gaussian matrix model with a source (\ref{17.1b}) as ME${}_{\beta,p}(x^{ a} e^{-x};\vec{m})$. In
 terms of this notation, the duality formula of \cite{De08} reads
 \begin{align}\label{Dr2}
  \Big \langle \prod_{j=1}^n \prod_{k=1}^p (s_j + {2 \over \beta} x_k) & \Big \rangle_{{\rm ME}_{\beta,p}(x^{a} e^{-x}; \vec{m})}
  \nonumber \\
  & = 
  \Big \langle \prod_{j=1}^n \prod_{k=1}^p (x_j  +  {2 \over \beta}  m_k) \Big  \rangle_{{\rm ME}_{4 /  \beta,n}(x^{(2/\beta)( a + 1) - 1} e^{-x}; \vec{s})}.
  \end{align}
 A verification of (\ref{Dr2}) analogous to that just given for (\ref{Dr1}) is possible, again distinct from the derivation given in \cite{De08}.
 Our working makes use of the generalized Laguerre polynomials $P_\kappa^{(L)}(\vec{x};\alpha;a)$ \cite{BF97a}, specified by the requirements that
 their expansion in terms of Jack polynomials has the structure
 $$
 P_\kappa^{(L)}(\vec{x};\alpha;a) = P_\kappa(\vec{x};\alpha) +
  \sum_{\mu: |\mu| < |\kappa|}
  \tilde{a}_{\kappa \mu} P_\mu(\vec{x};\alpha)
$$
for some coefficients $ \tilde{a}_{\kappa \mu}$, and that they exhibit the orthogonality
$$
    \Big \langle P_\kappa^{(L)}(\vec{x};\alpha;a) P_\rho^{(L)}(\vec{x};\alpha;a)  \Big \rangle_{{\rm ME}_{2/ \alpha,N}(x^ae^{-x})}
    \propto \delta_{\kappa,\rho}.
 $$
 With $\alpha = 2/\beta$, we seek to expand both sides of (\ref{Dr2}) as a series in $\{P_\mu^{(L)}(-\vec{s};1/\alpha;a+\alpha - 1) \}$.
 
 For this we require the Jack polynomial expansion formula (\ref{20a}), and the Laguerre analogue of the expansion formula (\ref{20}),
 \begin{equation}\label{23}
 \prod_{k=1}^N \prod_{l=1}^n (x_k - y_l) =
  \sum_{\kappa} (-1)^{|\mu|} P_\kappa^{(L)}(\vec{x};\alpha;a)
  \alpha^{-|\mu|/2} P_\mu^{(L)}(\sqrt{\alpha} \vec{y};1/\alpha;\alpha (a + 1) - 1),
  \end{equation}
which is derived in the same way as is indicated for (\ref{20}).  We require too the Laguerre analogues of the integration formulas
(\ref{20.1}) and (\ref{20.2}) \cite{BF97a},
 \begin{align}
 \label{23.1}
 \langle P_\kappa(\vec{x};\alpha) \rangle_{{\rm ME}_{2/\alpha,N}(x^a e^{-x}; \vec{z})} & =  (-1)^{|\kappa|} P_\kappa^{(L)}(-\vec{z};\alpha;a) \\
 \label{23.2}
  \langle P_\kappa^{(L)}(\vec{x};\alpha;a) \rangle_{{\rm ME}_{2/\alpha,N}(x^a e^{-x}; \vec{z})} & =  P_\kappa(\vec{z};\alpha).
  \end{align}
  
  Rewriting (\ref{23}) to read
  $$
  \prod_{j=1}^n \prod_{k=1}^p \Big (s_j + {2 \over \beta} x_k \Big ) =
  \sum_{\kappa} (-1)^{|\mu|} \Big ( {2 \over \beta} \Big )^{|\kappa|} P_\kappa^{(L)}(\vec{x};2/\beta;\beta a /2) 
  P_\kappa^{(L)}(-\vec{s};\beta/2;a + 2/\beta - 1),
  $$
  then substituting in the LHS of (\ref{Dr2}) and performing the integration over $\vec{x}$ using (\ref{23.2}) shows that the LHS of (\ref{Dr2})  can be rewritten
  \begin{equation}\label{d1}
  \sum_\kappa (-1)^{|\mu|} P_\kappa(2  \vec{m}/\beta;2/\beta)
  P_\mu^{(L)}(-\vec{s};\beta/2;a+2/\beta - 1).
  \end{equation}
  To express the RHS in an analogous series, we rewrite (\ref{20a}) to read
  $$
  \prod_{j=1}^n \prod_{k=1}^p ( x_j + 2 m_k/ \beta ) = \sum_\kappa P_\kappa(  2 \vec{m}/\beta ;2/\beta) 
  P_\mu(\vec{x};\beta/2).
  $$
  After performing the integration over $\vec{x}$ using (\ref{23.1}) we reproduce (\ref{d1}).
  
  \section{The averaged characteristic polynomial for Gaussian and chiral Gaussian ensembles with
a source: general $\beta > 0$}
\setcounter{equation}{0}
In this section the duality formulas (\ref{Dr1}) and (\ref{Dr2}) will be used to generalize (\ref{pp1}) and (\ref{pp1A}) to the
case of the Gaussian and chiral Gaussian ensembles with
a source for general $\beta > 0$. The resulting formulas will then be used to deduce scaled limits generalizing (\ref{2}).

\subsection{The Gaussian ensemble with a source}
Consider (\ref{Dr1}) with
$$
n = 1, \quad s_1 = \lambda, \quad 2/\alpha = \beta, \quad \vec{x} \mapsto i \vec{x}.
$$
Simple manipulation shows
\begin{equation}\label{fr}
\Big \langle \prod_{k=1}^N (\lambda - y_k) \Big \rangle_{{\rm ME}_{\beta,N}(e^{-\beta y^2/2};\vec{x})}
= \Big \langle \prod_{j=1}^N (\lambda - x_j + i y) \Big \rangle_{{\rm ME}_{y \in {\rm N}[0,1/\sqrt{2}]}},
\end{equation}
where ME${}_{\beta,N}(e^{-\beta y^2/2};\vec{x})$ refers to the PDF (\ref{AF1}) with $\tilde{\beta} = \beta/2$, $t=1$, $\vec{\lambda} = \vec{y}$,
$\vec{\mu}=\vec{x}$. For $\beta = 1$ ($\beta = 2$) this
is precisely (\ref{pp1}) in the case that $H$ has real (complex) off-diagonal elements.

The RHS of (\ref{fr}) is independent of $\beta$. But in the case $\beta = 2$ and $\vec{x} \mapsto \vec{x}/2$ we know that the average
characteristic polynomial gives the multiple Hermite polynomials of type II \cite{BK03}. Thus this identification persists for general $\beta > 0$.

In relation to the asymptotics we restrict attention to the case $x_j = 0$ $(j=r+1,\dots,N)$. Minor manipulation of (\ref{fr}) shows
$$
\Big \langle \prod_{k=1}^N (\lambda - y_k) \Big \rangle_{{\rm ME}_{\beta,N}(e^{-\beta y^2/2};\vec{x})} = (-1)^N {\sqrt{\pi} \over 2^N} e^{\lambda^2}
\Gamma^{(r+1)}(\lambda;\{-2x_k\}) ,
$$
where
\begin{equation}\label{Hp}
\Gamma^{(r+1)}(u;\{a_k\}) :=
\int_{-i \infty}^{i \infty} e^{y^2/4 + u y }y^{N-r} \prod_{j=1}^r ( y - a_j) \, {dy \over 2 \pi i}.
\end{equation}
The function $\Gamma^{(r+1)}(u;\{a_k\})$ is a particular incomplete multiple Hermite polynomial \cite{DF06b}.

Introducing the scaled variables
\begin{equation}\label{scaled}
\lambda = \sqrt{2N} + {X \over \sqrt{2} N^{1/6}}, \qquad
x_k = \sqrt{N/2} - {N^{1/6} s_k \over \sqrt{2}} \quad (k=1,\dots,r)
\end{equation}
the leading large $N$ asymptotic form of the  incomplete multiple Hermite polynomial can be read off from
\cite[Proof of Prop.~18]{DF06b}. This involves the incomplete multiple Airy function
\begin{align}\label{Ai}
(-1)^{r+1} {\rm Ai}^{(r+1)}(X,\{s_k\}) & := \int_{\mathcal A} e^{- X w + w^3/3}
\prod_{k=1}^r (w + s_k) \, {dw \over 2 \pi i} \nonumber \\
& = \prod_{k=1}^r \Big ( - {\partial \over \partial X} + s_k \Big ) {\rm Ai}(X),
\end{align}
where $\mathcal A$ is a simple contour, starting at $e^{-\pi i /3} \infty$, and finishing at $e^{\pi i /3} \infty$, following the
corresponding rays asymptotically. As a consequence, we can generalize (\ref{2}) to the case of the Gaussian ensemble with a source.

\begin{proposition}
Let $\lambda$ and $x_k$ be given in terms of $X$ and $s_k$ according to (\ref{scaled}). Define
$$
C_N^{(r+1)} = 2^{-(N-1)/2}N^{(N+1)/2 - (r+1)/3} e^{-N/2}.
$$
We have
\begin{equation}\label{fna}
\lim_{N \to \infty} {e^{-\lambda^2/2} \over C_N^{(r+1)}} \Big \langle \prod_{k=1}^N (\lambda - y_k) \Big \rangle_{{\rm ME}_{\beta,N}(e^{-\beta y^2/2};\vec{x})} =
(-1)^{r+1} {\rm Ai}^{(r+1)}(X,\{s_k\}) .
\end{equation}
\end{proposition}

Making use of Stirling's formula to replace $(N!)^{1/2}$ in (\ref{2}), and noting from (\ref{Ai}) that $- {\rm Ai}^{(1)}(X,\{s_k\})  =  {\rm Ai}(X)$ shows
that the $r=0$ case of (\ref{fna}) agrees with (\ref{2}).

\subsection{The chiral Gaussian ensemble with a source}
With $n=1$, $s_1 = - \lambda$, $\vec{m} \mapsto \beta \vec{m}/2$, $a \mapsto (\beta/2)(a+1) - 1$ and the change of variables
$2 x_k/ \beta \mapsto x_k$ on the LHS, the duality formula (\ref{Dr2}) reads
\begin{eqnarray}\label{box}
&&\Big \langle \prod_{k=1}^p (\lambda - x_k) \Big \rangle_{{\rm ME}_{\beta,p}(x^{(\beta/2)(a+1)-1} e^{-\beta x/2};\beta \vec{m}/2)}  \nonumber \\
&& \qquad =
{(-1)^p e^\lambda \over \Gamma(a+1)}
\int_0^\infty x^a e^{-x} \, {}_0 F_1 (a+1; - \lambda x) \prod_{k=1}^p (x + m_k) \, dx.
\end{eqnarray}
We recognize the RHS of (\ref{box}) as the multiple Laguerre polynomial of type II \cite[eq.~(77)]{DF07}.
Also, we remark that
this evaluation of the averaged characteristic polynomial, for $a = n - p$ and $\beta = 1$ or 2 is in precise agreement with
the result of Corollary \ref{C3}.

A systematic study of the asymptotics of the multiple Laguerre polynomials is yet to be undertaken. We will initiate this
task in the simplest nontrivial case of 
\begin{equation}\label{nt}
m_2 = \cdots = m_p = 0.
\end{equation}

\begin{proposition}
Let (\ref{nt}) hold, and introduce $X$ and $s_1$ so that
\begin{equation}\label{Xm}
\lambda = 4p + 2a + 2 + 2(2p)^{1/3} X, \quad m_1 = p - (2p)^{2/3} s_1.
\end{equation}
Also define $D_p^{(2)} = p! (2p)^{1/3} 2^{-a}$. We have
\begin{equation}\label{DD}
\lim_{p \to \infty} {e^{-\lambda/2} \over D_p^{(2)}}  \Big \langle \prod_{k=1}^p (\lambda - x_k) \Big \rangle_{{\rm ME}_{\beta,p}(x^{(\beta/2)(a+1)-1} e^{-\beta x/2};\beta \vec{m}/2)}
= s_1 {\rm Ai}(X) - {\rm Ai}'(X).
\end{equation}
\end{proposition}

\noindent
{\it Proof.} 
Recalling (\ref{LP}) we see that the RHS of (\ref{box})
can be written in terms of Laguerre polynomials to read
\begin{equation}\label{DD1}
(-1)^p p! L_p^a(\lambda) - m_1 (-1)^{p-1} (p-1)! L_{p-1}^a(\lambda).
\end{equation}
 But it is a consequence of an asymptotic formula in Szeg\"o's book \cite[p.~201]{Sz75} that for fixed $k$
 and with $\lambda$ as in (\ref{Xm})
 $$
 e^{-\lambda/2} (-1)^{p+k} L_{p+k}^a(\lambda) =
 2^{-a - 1/3} p^{-1/3} \Big ( {\rm Ai}(X) - {2k \over (2p)^{1/3}} {\rm Ai}'(X) + {\rm O}(p^{-2/3}) \Big ).
 $$
 Multiplying (\ref{DD1}) by $e^{-\lambda/2}$, substituting for $m_1$ according to (\ref{Xm}) and using this
 asymptotic formula with $k=0$ and $k=-1$ gives (\ref{DD}).
 \hfill$\square$
 
 \medskip
 Note from (\ref{Ai}) that the scaled limit (\ref{DD}) is the same as the scaled limit (\ref{fna}) in the case $r=1$ --- an example
 of universality. More generally, choosing $m_k = p - (2p)^{2/3} s_k$, $k=1,\dots,r$, and $m_k =0 $, $k=r+1,\dots,p$ we would expect that for
 an appropriate $D_p^{(r+1)}$, that upon multiplication by $ e^{-\lambda/2} / D_p^{(r+1)} $, the large $p$ limit of the
 averaged characteristic polynomial in (\ref{box}) will equal the RHS of (\ref{fna}). 
 
 \subsection{Related work and future problems}
 In the Introduction it was remarked that, in the classical setting,  one link with integrable systems is provided by the fact the Hermite polynomial
 in (\ref{1}) satisfies a differential equation. It is also true that the incomplete Hermite polynomial (\ref{Hp}) satisfies a linear differential equation ---
 now of degree $r+2$ (see e.g.~\cite[eq.~(23.8.7)]{Is05}). A future challenge is to find a derivation of this differential equation, working directly
 from the definition of the random matrix average in (\ref{fr}). Such a derivation is quite straightforward in the case of the average in (\ref{1}) \cite{Ao87},
 upon making use of the explicit form of the joint eigenvalue PDF, and is in fact perhaps the simplest result in the broader concern of finding (partial) differential equations associated with
 Selberg correlation functions \cite{Ka93}, \cite{Fo11A}.
 
 Another way to characterize the Hermite polynomial in (\ref{1}) is via its three term recurrence. The tridiagonal construction of the Gaussian
 $\beta$-ensemble due to Dumitriu and Edelman leads most naturally to a direct demonstration that the characteristic polynomial on the
 LHS of (\ref{1}) satisfies the same three term recurrence \cite{DE02}. It is similarly the case that the multiple classical polynomials corresponding
 to the RHSs of (\ref{fr}) and (\ref{box}) satisfy linear recurrences \cite{vA11}, thus suggesting the problem of a direct derivation of the same
 recurrences for the  characteristic polynomials on the LHSs. For the shifted mean Gaussian ensembles of Section \ref{s2.1}, with
 $s_{r+1},\dots, s_N = 0$, the $r  \times r$ block tridiagonal form found recently by Bloemendal and Vir\'ag \cite{BV11}, may provide some
 insight. 
 
 The duality formulas (\ref{Dr1}) and (\ref{Dr2}) give formulas well suited to asymptotic analysis, not just for the averaged
 characteristic polynomial, but for the average of an arbitrary (but fixed) number of characteristic polynomials and thus by
 coalescence for integer moments of the characteristic polynomials. Using the duality formula (\ref{w2}) the
 soft edge scaling of the latter for the Gaussian $\beta$-ensemble was given in terms of a $\beta$-generalization of the
 Kontsevich matrix Airy function \cite{DF06}, and for $\beta = 1,2$ and 4 (general $\beta > 0$) this was generalized in \cite{De08} (\cite{DL11})
 to the case with the variables
 not coalesced. Of course one would like to compute the generalization of the soft edge scaling limit (\ref{fna}) from a single characteristic
 polynomial to a product of characteristic polynomials.
 In has been announced by Desrosiers and Liu in \cite{DL11} that they are preparing a paper tackling this problem.

 \subsection*{Acknowledgement}
 This work was supported by the Australian Research Council.
 
 \appendix
 \section{A recursive construction of the Gaussian $\beta$-ensemble with a source}
\setcounter{equation}{0}
The Gaussian and chiral Gaussian matrix models with a source are specified by the PDFs proportional to (\ref{AF1}) and
(\ref{17.1b}) respectively. We have seen that both can be interpreted for general $\beta > 0$ as Green function solutions of
particular Fokker-Planck equations describing the Brownian evolution of the underlying log-gases. Here we will specify the recursive
construction of a random matrix ensemble having the eigenvalue PDF (\ref{AF1}). While no such construction is known in relation to
(\ref{17.1b}), the PDF specifying the general variance Laguerre $\beta$-ensemble can similarly be generated as the eigenvalue PDF
of recursively defined random matrices \cite{Fo11y}.

Consider random matrices specified by the distribution on elements (\ref{8.0}), but with the $\tau$-dependence scaled out so that it
reads
\begin{equation}\label{8.a}
e^{-(\beta/2) {\rm Tr} (H - M)^2}.
\end{equation}
Let us take $\beta = 1$ (real entries) for definiteness, and furthermore write $Y = Y_N$ to indicate that the matrix $Y$ has size
$N \times N$. Then by first noting that (\ref{8.a}) can equivalently be written
\begin{equation}
e^{-(\beta/2) {\rm Tr} (H_N - D_N)^2},
\end{equation}
where $D_N = {\rm diag} \, M_N = {\rm diag} \, (\mu_1,\dots,\mu_N)$, as observed recently in \cite{BFF10} we have that
\begin{equation}\label{8.b}
{\rm ePDF} \, H_{N+1} =
{\rm ePDF} \begin{bmatrix}
x_{11} & x_{12} & \cdots & x_{1 \, N+1} \\
x_{12} &  \lambda_1^{(1)} & \cdots & 0 \\
\vdots & & \ddots & \vdots \\
x_{1 \, N+1} & 0 & \hdots & \lambda_N^{(N)}
\end{bmatrix}.
\end{equation}
Here ePDF refers to the eigenvalue PDF, $\{\lambda_j^{(N)} \}_{j=1,\dots,N}$ denotes the eigenvalues of $H_N$ and
$$
x_{11} \mathop{=}\limits^{\rm d} {\rm N}[\mu_{N+1},1], \qquad
x_{1 \, j}  \mathop{=}\limits^{\rm d} {\rm N}[0,1/\sqrt{2}] \quad (j=2,\dots,N+1).
$$
In fact in (\ref{8.b}) the eigenvalue PDF of $H_{N+1}$ is fully determined by the distribution
\begin{equation}\label{8.c}
| x_{1j}|^2   \mathop{=}\limits^{\rm d}  \Gamma[\beta/2;1],
\end{equation}
where $\Gamma[s;\sigma]$ refers to the gamma distribution (PDF proportional to $x^{s-1} e^{-s/\sigma}$, $s > 0$).
This remains true for $\beta = 2$ and 4 (the RHS of (\ref{8.b}) should then be scaled by $1/\sqrt{\beta}$, but this detail can be ignored).

But (\ref{8.c}) makes sense for general $\beta > 0$, and allows us to use (\ref{8.b}) to give a recursive generation of the eigenvalue PDF
$\beta$-ensemble coinciding with the eigenvalue PDF of the shifted mean Gaussian ensembles for $\beta = 1,2$ and 4. The use of (\ref{8.c})
to $\beta$-generalize classical random matrix ensembles goes back to \cite{DE02}), while the idea of a recursive construction can be found in
\cite{FR02b}. According to (\ref{AF1}) the eigenvalue PDF of $H_N$ is proportional to
\begin{equation}\label{AF2}
\prod_{1 \le j < k \le N} |\lambda_k - \lambda_j |^\beta
e^{-\sum_{j=1}^N \lambda_j^2/2 - \sum_{j=1}^N \mu_j^2/2}
 {}^{\phantom{/}}_0 {\mathcal F}_0^{(2/\beta)} (\vec{\lambda};  \vec{\mu})
 \end{equation}
 for $\beta = 1,2$ and 4. We would like to show that the recursive construction (\ref{8.b}) implies that (\ref{AF2}) remains true for
 general $\beta > 0$. Our strategy will be to show that (\ref{8.b}) implies a recurrence formula in $N$ for the eigenvalue PDF, and that
 (\ref{AF2}) satisfies the recurrence.
 
 \begin{proposition}
 Consider the random matrix on the RHS of (\ref{8.b}). Let the eigenvalues be denoted $\{\lambda_j^{(N+1)}\}_{j=1,\dots,N+1}$.
 With $\{\lambda_j^{(N)}\}_{j=1,\dots,N}$ given, the eigenvalue PDF is proportional to
 \begin{align}\label{2H}
 & e^{-(\sum_{j=1}^{N+1} (\lambda_j^{(N+1)})^2 - \sum_{j=1}^N (\lambda_j^{(N)})^2)}
e^{ \mu_{N+1} (\sum_{j=1}^{N+1} \lambda_j^{(N+1)} - \sum_{j=1}^N  \lambda_j^{(N)})} \nonumber \\
& \quad \times
{\prod_{1 \le j < k \le N + 1} (\lambda_j^{(N+1)} - \lambda_k^{(N+1)} ) \over
\prod_{1 \le j < k \le N } (\lambda_j^{(N)} - \lambda_k^{(N)} )^{\beta - 1} }
\prod_{j=1}^N \prod_{k=1}^{N+1} | \lambda_j^{(N)} - \lambda_k^{(N+1)} |^{\beta/2 - 1},
\end{align}
subject to the interlacing 
\begin{equation}\label{2.7}
\lambda_1^{(N+1)} > \lambda_1^{(N)} > \lambda_2^{(N+1)} > \lambda_2^{(N)} > \cdots > \lambda_N^{(N)} > \lambda_{N+1}^{(N+1)}.
\end{equation}
\end{proposition}

\noindent
{\it Proof.} \quad It is a simple exercise to show that the secular equation for the matrix on the RHS of (\ref{8.b}) is
\begin{equation}\label{2.6}
0 = \lambda - x_{11} - \sum_{j=1}^N {q_j \over \lambda - \lambda_j^{(0)}}, \qquad q_j = | x_{1,j+1} |^2.
\end{equation}
The distribution of the zeros of this random rational function in the case of $x_{11}$ fixed is given by \cite{FR02b}
(Prop.~8 with $n = N$, $s_1 = \cdots = s_N = \beta/2$, $\lambda \mapsto \lambda - x_{11}$, $a_j \mapsto a_j - x_{11}$).
Explicitly we have that in this setting the PDF for the zeros is supported on the region (\ref{2.7}) and is subject to the further constraint
\begin{equation}\label{2.7a}
\sum_{l=1}^{N+1} \lambda_l^{(N+1)} = \sum_{l=1}^N \lambda_l^{(N)} + x_{11}.
\end{equation}
Subject to these requirements, its functional form is proportional to 
 \begin{align}
 & 
e^{-(\sum_{j=1}^{N+1} (\lambda_j^{(N+1)})^2 - \sum_{j=1}^N (\lambda_j^{(N)})^2)} e^{x_{11}^2/2} \nonumber \\
& \quad
\times 
{\prod_{1 \le j < k \le N + 1} (\lambda_j^{(N+1)} - \lambda_k^{(N+1)} ) \over
\prod_{1 \le j < k \le N + 1} (\lambda_j^{(N)} - \lambda_k^{(N)} )^{\beta - 1} }
\prod_{j=1}^N \prod_{k=1}^{N+1} | \lambda_j^{(N)} - \lambda_k^{(N+1)} |^{\beta/2 - 1}.
\end{align}
Multiplying this by the PDF of $x_{11}$, expanding the exponent and making use of (\ref{2.7a}) gives (\ref{2H}).
\hfill $\square$

\medskip
\begin{corollary}
The joint PDF of $\{\lambda_j^{(N)}\}_{j=1,\dots,N}$ and $\{\lambda_j^{(N+1)}\}_{j=1,\dots,N}$ is proportional to
\begin{align}\label{Sat}
 & 
e^{-(\sum_{j=1}^{N+1} (\lambda_j^{(N+1)})^2/2}
e^{\mu_{N+1}( \sum_{j=1}^{N+1} \lambda_j^{(N+1)} - \sum_{j=1}^N \lambda_j^{(N)})} \prod_{1 \le j < k \le N + 1}
(\lambda_j^{(N+1)} - \lambda_k^{(N+1)}) \nonumber \\
& \quad \times
\prod_{1 \le j < k \le N} (\lambda_j^{(N)} - \lambda_k^{(N)} ) \prod_{j=1}^N \prod_{k=1}^{N+1}
| \lambda_j^{(N)} - \lambda_k^{(N+1)}|^{\beta/2 - 1} 
 {}^{\phantom{/}}_0 {\mathcal F}_0^{(2/\beta)} (\vec{\lambda}^{(N)};  \vec{\mu}^{(N)})
\end{align}
subject to the interlacing (\ref{2.7}).
\end{corollary}

\noindent
{\it Proof.} \quad We simply multiply together (\ref{2H}) and the PDF ME${}_{\beta,N}(e^{-x^2})$.
\hfill $\square$

\medskip
To proceed further requires knowledge of (a special case of) the so-called Dixon-Anderson density
(see e.g.~\cite[eq.~(4.11)]{Fo10}). This is the conditional PDF for $\{\lambda_j^{(N)} \}$ given
$\{\lambda_j^{(N+1)} \}$ (note that (\ref{2H}) is the other way around) specified by
\begin{equation}\label{DA}
{\Gamma((N+1)\beta/2) \over (\Gamma(\beta/2))^{N+1}}
{\prod_{1 \le j < k \le N } (\lambda_j^{(N)} - \lambda_k^{(N)} ) \over
\prod_{1 \le j < k \le N + 1} (\lambda_j^{(N+1)} - \lambda_k^{(N+1)} ) }
\prod_{j=1}^N \prod_{k=1}^{N+1} | \lambda_j^{(N)} - \lambda_k^{(N+1)} |^{\beta/2 - 1},
\end{equation}
and again subject to the interlacing (\ref{2.7}). Referring to this as DA${}_{N+1}(\beta/2)$, it is a known
result that \cite{GK02} (see also \cite[Appendix C]{FR02b})
\begin{eqnarray}\label{16.1}
&&
e^{\mu_{N+1} \sum_{j=1}^{N+1}  \lambda_j^{(N+1)}}
\langle e^{- \mu_{N+1} \sum_{j=1}^{N} \lambda_j^{(N)}}
 {}_0^{} {\mathcal F}_0^{(2/\beta)}(\lambda^{(N)};\mu^{(N)})
 \rangle_{{\rm DA}_N(2/\beta)} \nonumber \\&& \qquad =
 {}_0^{} {\mathcal F}_0^{(2/\beta)}(\vec{\lambda}^{(N+1)};\mu^{(N+1)}). 
 \end{eqnarray}

On the other hand, it follows from (\ref{Sat}) that the $N \mapsto N + 1$ version of (\ref{AF2}) is also the PDF of
$\{\lambda_j^{(N+1)}\}$ for general $\beta > 0$, if we can show that
\begin{eqnarray}\label{15.1}
 && \int d \vec{\lambda}^{(N)} \,
e^{\mu_{N+1} ( \sum_{j=1}^{N+1} \lambda_j^{(N+1)} - \sum_{j=1}^N \lambda_j^{(N)})}
\prod_{1 \le j < k \le N} (\lambda_j^{(N)} - \lambda_k^{(N)} )
\prod_{j=1}^N \prod_{k=1}^{N+1} | \lambda_j^{(N)} - \lambda_k^{(N+1)}|^{\beta/2 - 1} \nonumber \\
 && \quad \times  {}_0^{} {\mathcal F}_0^{(2/\beta)}(\vec{\lambda}^{(N+1)};\mu^{(N+1)})  \nonumber \\
 && \qquad
\propto \prod_{1 \le j < k \le N + 1} ( \lambda_j^{(N+1)} -  \lambda_k^{(N+1)} )^{\beta - 1}
  {}_0^{} {\mathcal F}_0^{(2/\beta)}(\vec{\lambda}^{(N+1)};\mu^{(N+1)}). 
\end{eqnarray}
But (\ref{16.1}) is precisely (\ref{15.1}), complete with normalization. It thus follows that (\ref{AF2}) is
indeed the eigenvalue PDF of $H_{N+1}$.


\providecommand{\bysame}{\leavevmode\hbox to3em{\hrulefill}\thinspace}
\providecommand{\MR}{\relax\ifhmode\unskip\space\fi MR }
\providecommand{\MRhref}[2]{%
  \href{http://www.ams.org/mathscinet-getitem?mr=#1}{#2}
}
\providecommand{\href}[2]{#2}

\end{document}